\documentclass[10pt]{article}
\usepackage{amsmath,amssymb,mathdots, graphicx, color}
\usepackage[margin=1in]{geometry}
\usepackage[affil-it]{authblk}

\newmuskip\pFqskip
\mathchardef\pFcomma=\mathcode`, 

\newcommand*\pGq[5]{%
  \begingroup
  \begingroup\lccode`~=`,
    \lowercase{\endgroup\def~}{\pFcomma\mkern\pFqskip}%
  \mathcode`,=\string"8000
  {}_{#1}\phi_{#2}\biggl[\genfrac..{0pt}{}{#3}{#4}\bigg|#5\biggr]%
  \endgroup
}

\numberwithin{equation}{section}

\title{A $q$-generalization of the para-Racah polynomials}
\author[1]{Jean-Michel Lemay}
\affil[1]{Centre de Recherches Math\'ematiques, Universit\'e de Montr\'eal, C.P. 6128,\protect\\ Succ. Centre-ville, Montr\'eal, QC, Canada, H3C 3J7}
\author[1]{Luc Vinet}
\author[2]{Alexei Zhedanov}
\affil[2]{Department of Mathematics, School of Information, Renmin University of China, Beijing 100872, China}
\date{}

\begin{document}
\maketitle

\begin{abstract}
New bispectral orthogonal polynomials are obtained from an unconventional truncation of the Askey-Wilson polynomials. In the limit $q \to 1$, they reduce to the para-Racah polynomials which are orthogonal with respect to a quadratic bi-lattice. The three term recurrence relation and q-difference equation are obtained through limits of those of the Askey-Wilson polynomials. An explicit expression in terms of hypergeometric series and the orthogonality relation are provided. A $q$-generalization of the para-Krawtchouk polynomials is obtained as a special case. Connections with the $q$-Racah and dual-Hahn polynomials are also presented.
\end{abstract}

\section{Introduction}

Persymmetric Jacobi matrices are invariant under reflections with respect to the anti-diagonal \cite{DEBOOR1978245,0266-5611-3-4-010,de1986finite,2004_Gladwell}. Recent studies found applications for these matrices in transfer of quantum information along Heisenberg spin chains \cite{2005_Shi&Li&Song&Sun_PhysRevA_71_032309,2010_Chakrabarti&VanderJeugt_JPhysA_43_085302,2010_Jafarov&VanderJeugt_JPhysA_43_405301,2011_Stoilova_Van_der_Jeugt_Hahn,2012_Vinet&Zhedanov_PhysRevA_85_012323,2012_Vinet&Zhedanov_JPhysA_45_265304,2012_Vinet&Zhedanov_JPhysConfSer_343_012125}. Necessary and sufficient conditions for achieving perfect end-to-end transfer of a qubit have been framed as properties of Jacobi matrices and their spectra \cite{2012_Vinet&Zhedanov_PhysRevA_85_012323}. In addition to persymmetry, differences of consecutive eigenvalues must satisfy certain relations. This suggested to study Jacobi matrices whose spectrum are the superposition of two lattices (bi-lattices). Interestingly, this idea led to the characterization of two novel sets of orthogonal polynomials : the para-Krawtchouk polynomials \cite{2012_Vinet&Zhedanov_JPhysA_45_265304} which are orthogonal with respect to a linear bi-lattice and the para-Racah polynomials \cite{lemay2016analytic} which are orthogonal with respect to a quadratic bi-lattice.

It turns out that the para-Racah polynomials can be seen to arise from a singular truncation of the Wilson polynomials which sits atop the Askey scheme \cite{lemay2016racah}. They are however not classical in the usual sense since they satisfy the same difference equation as the Wilson polynomials but with degenerate eigenvalues. In this setting, the corresponding Jacobi matrix is not persymmetric, but rather corresponds to an one-parameter isospectral deformation of the persymmetric one. Hence, the general para-Racah polynomials depend on three real parameters $a,c,\alpha$ and the Jacobi matrix is persymmetric only when $\alpha=1/2$. Additionally, the para-Krawtchouk polynomials can be recovered from the para-Racah polynomials by an appropriate limit. 

Our goal here is to report on a $q$-generalization of these results. More precisely, we will construct new orthogonal polynomials via a singular truncation of the Askey-Wilson polynomials. We name them the $q$-para-Racah polynomials since they reduce to the para-Racah ones in the limit $q\to1$. A study of the most general tridiagonal representations of the $q$-oscillator algebra $AB-qBA=1$ has given hints of their existence and their connection to the $q$-para-Krawtchouk polynomials \cite{tsujimoto2017tridiagonal}. 

The paper will unfold as follows. In section 2, we review some basic properties of the Askey-Wilson polynomials. We then proceed with the construction of a set of $N+1$ $q$-para-Racah polynomials. We start in section 3 with the case of $N=2j+1$. The three-term recurrence relation, the $q$-difference equation and an explicit expression in terms of hypergeometric series are obtained from a singular truncation of the Askey-Wilson OPs. The general orthogonality relation is obtained by making use of persymmetry features that are observed when $\alpha=1/2$. Section 4 is dedicated to the case $N=2j$ and the corresponding formulas. A $q$-generalization of the para-Krawtchouk polynomials is obtained as a special case in section 5. A connection to the $q$-Racah and dual-Hahn polynomials is also presented and a short conclusion follows.

\section{Askey-Wilson polynomials}

  Let us first review some properties of the Askey-Wilson polynomials. They are sitting atop the continuous part of the $q$-Askey scheme and depend on four parameters $a,b,c,d$. We shall denote them by $W_n(x;a,b,c,d|q)$ or simply $W_n(x)$ when the parameters need not be explicit. They admit a simple expression in terms of the hypergeometric function  
  \begin{align}
  \begin{aligned} \label{AWPolynomials}
  W_n(x;a,b,c,d|q)=\pGq{4}{3}{q^{-n},abcdq^{n-1},ae^{i\theta},ae^{-i\theta}}{ab,ac,ad}{q;q}
  \end{aligned}
  \end{align}
  in the variable $x = \cos(\theta)$. For future considerations, it will be useful to write this hypergeometric series explicitly 
  \begin{align} \label{AWExplicit}
  W_n(x;a,b,c,d|q) \equiv \sum_{k=0}^\infty A_{n,k}\Phi_k(x)
  \end{align}
  where
  \begin{align}
  \begin{aligned} \label{AWExplicitCoefficients}
  A_{n,k} = \frac{(q^{-n},abcdq^{n-1};q)_k}{(q,ab,ac,ad;q)_k},\quad \Phi_k(x) = (ae^{i\theta},ae^{-i\theta};q)_k,
  \end{aligned}
  \end{align} 
  with the standard notation for the $q$-Pochhammer symbol $(a;q)_k = (1-a)(1-aq)\dots(1-aq^{k-1})$ and its products $(a_1,a_2,\dots,a_n;q)_k=(a_1;q)_k (a_2;q)_k\dots(a_n;q)_k$.
  The Askey-Wilson polynomials satisfy a $q$-difference equation 
  \begin{align} \label{AWdifferenceeqn}
   q^{-n}(1-q^n)(1-abcdq^{n-1})W_n(x) = A(\theta)T_+W_n(x) -[A(\theta)+\overline{A}(\theta)]W_n(x) + \overline{A}(\theta)T_-W_n(x)
  \end{align}
  where the $T_\pm$ are the following $q$-shift operators
  \begin{align} \label{qshiftop}
  \begin{aligned}
   &T_+e^{i\theta} = qe^{i\theta} & &T_-e^{i\theta} = q^{-1}e^{i\theta} \\
   &T_+e^{-i\theta} = q^{-1}e^{-i\theta}\quad & &T_-e^{-i\theta} = qe^{-i\theta} 
  \end{aligned}
  \end{align}
  and the coefficients $A(\theta)$ and its complex conjugate $\overline{A}(\theta)$ are defined by
  \begin{align}
   A(\theta) = \frac{(1-a e^{i\theta})(1-b e^{i\theta})(1-c e^{i\theta})(1-d e^{i\theta})}{(1- e^{2i\theta})(1-q e^{2i\theta})}.
  \end{align}
  Moreover, the $W_n(x)$ obey the three-term recurrence relation 
  \begin{align}
  2xW_n(x)=A_nW_{n+1}(x)+(a+a^{-1}-A_n-C_n)W_n(x)+C_nW_{n-1}(x)
  \end{align}
  with recurrence coefficients
  \begin{align} \label{rec_coef_wilson}
  \begin{aligned}
  A_n =& \frac{(1-abq^n)(1-acq^n)(1-adq^n)(1-abcdq^{n-1})}{a(1-abcdq^{2n-1})(1-abcdq^{2n})}, \\
  C_n =& \frac{a(1-q^n)(1-bcq^{n-1})(1-bdq^{n-1})(1-cdq^{n-1})}{(1-abcdq^{2n-2})(1-abcdq^{2n-1})}.
  \end{aligned}
  \end{align}
  It is sometimes more convenient to work with the monic Askey-Wilson polynomials $\tilde{W}_n(x)$. They are given by 
  \begin{align}
   W_n(x) = \frac{2^n}{A_1A_2\dots A_{n-1}} \tilde{W}_n(x)
  \end{align}
  and satisfy the recurrence relation
  \begin{align} \label{AWrecrel}
  x\tilde{W}_n(x) = \tilde{W}_{n+1}(x) + \tfrac12(a+a^{-1}-A_n-C_n)\tilde{W}_{n}(x) + \tfrac14 A_{n-1}C_n\tilde{W}_{n-1}(x).
  \end{align}
  The orthogonality relation and many more properties on the Askey-Wilson polynomials can be found in \cite{2010_Koekoek_&Lesky&Swarttouw,Ismail-2009}. It is well-known that the Askey-Wilson polynomials can be reduced to a finite set of $N+1$ orthogonal polynomials if
  \begin{align} \label{truncation_cond}
  A_{N}C_{N+1} = 0.
  \end{align}
  It is easy to see that there are multiple ways to achieve this condition by looking at the recurrence coefficients \eqref{rec_coef_wilson}. Many of these possibilities lead to the $q$-Racah polynomials. More precisely, if one chooses the parameters so that either $(1-abq^N),(1-acq^N),(1-adq^N),(1-bcq^{N}),(1-bdq^{N})$ or $(1-cdq^{N})$ is equal to zero, then the corresponding polynomials will be the $q$-Racah polynomials. There is however another possibility which consists in choosing 
  \begin{align} \label{truncation}
  1-abcdq^{N-1}=0,
  \end{align}
  but this introduces a singularity in the denominators of $A_n$ and $C_n$ for $n\sim N/2$. The present article deals with this option and we show that under an appropriate parametrization, it is possible to recover a set of orthogonal polynomials distinct from the $q$-Racah polynomials that has not been studied before as far as we know. In \cite{lemay2016racah}, we presented a similar argument for the case $q=1$.

\section{Odd case : $N = 2j+1$}  

We aim to study the singular truncation \eqref{truncation} and to characterize the corresponding orthogonal polynomials. Remark first that depending on the parity of $N$, the singular factors in the denominators of $A_n$ and $C_n$ given by \eqref{rec_coef_wilson} change. This implies that both cases must be treated separately. This section will deal with the case $N$ odd and we relegate the case $N$ even to the next section. 

Now, fix $N=2j+1$ and let 
\begin{align} \label{trunclimit}
 b = a^{-1}q^{-j+e_1t}, \qquad d = c^{-1}q^{-j+e_2t}, \qquad t\to0.
\end{align}
Notice that in the limit $t\to0$, this parametrization realizes the truncation condition \eqref{truncation}. 
\subsection{Recurrence relation}
Injecting the formulas \eqref{trunclimit} in the recurrence coefficients \eqref{rec_coef_wilson} and taking the limit $t\to0$, one can check that the only coefficients depending on $e_1$ and $e_2$ are
\begin{align}
 \lim_{t\to0} A_j = \left(\frac{e_1}{e_1+e_2}\right)\dfrac{(1-acq^j)(c-a)(1-q^{-j-1})}{ac(1-q^{-1})}, \\
 \lim_{t\to0} C_{j+1} = \left(\frac{e_2}{e_1+e_2}\right)\dfrac{(1-q^{j+1})(c-a)(ac-q^{-j})}{ac(1-q)}.
\end{align}
Remark that the two occurences of the parameters $e_1$ and $e_2$ are not independent since they sum to one. They can hence be combined into a single deformation parameter $\alpha$ :
\begin{align}
 \frac{e_1}{e_1+e_2} = \alpha, \qquad \frac{e_2}{e_1+e_2} = 1-\alpha.
\end{align}
With this notation, the limiting recurrence coefficients can be expressed as 
\begin{align} 
 \lim_{t\to0} A_n &= \begin{cases}
       \dfrac{(1-acq^n)(c-aq^{n-j})(1-q^{n-2j-1})}{ac(1-q^{2n-2j-1})(1+q^{n-j})} \qquad &n\neq j, \\[1em]
       \dfrac{\alpha(1-acq^j)(c-a)(1-q^{-j-1})}{ac(1-q^{-1})} \qquad &n=j,
       \end{cases}\label{reccoefA}\\[1em]
 \lim_{t\to0} C_n &= \begin{cases}
       \dfrac{(1-q^n)(a-cq^{n-j-1})(ac-q^{n-2j-1})}{ac(1+q^{n-j-1})(1-q^{2n-2j-1})} \qquad &n\neq j+1, \\[1em]
       \dfrac{(1-\alpha)(1-q^{j+1})(a-c)(ac-q^{-j})}{ac(1-q)} \qquad &n= j+1.
       \end{cases}\label{reccoefC}
\end{align}
These new recurrence coefficients verify $A_N C_{N+1}=0$ with $N=2j+1$ and thus provides a finite set of polynomials. We define the (monic) $q$-para-Racah polynomials $R_n(x;a,c,\alpha|q)$ (or $R_n(x)$ for short) via the three-term recurrence relation
\begin{align} \label{recrel1}
 &x R_n(x) = R_{n+1}(x) + b_nR_n(x) + u_n R_{n-1}(x) 
\end{align}
with initial conditions $R_{-1}(x)=0, R_0(x)=1$ and where the recurrence coefficients are those of the monic Askey-Wilson polynomials \eqref{AWrecrel} under the limit given by \eqref{trunclimit}
\begin{align}
 &b_n = \lim_{t\to0}\tfrac{1}{2}(a+a^{-1}-A_n-C_n), \\
 &u_n = \lim_{t\to0}\tfrac{1}{4}A_{n-1}C_n. 
\end{align}
Using formulas \eqref{reccoefA} and \eqref{reccoefC}, they can be written as
{\small
\begin{align}
 b_n &= \begin{cases}
        \dfrac{(a+c) \left(q^{j+1}+1\right) q^n \left(a c q^j+1\right)}{2 a c \left(q^j+q^n\right) \left(q^{j+1}+q^n\right)} \qquad &n \neq j,j+1, \\[1em]
        \dfrac{a+a^{-1}}{2}+\dfrac{\alpha(c-a)\left(q^{j+1}-1\right)q^{-j}\left(acq^j-1\right)}{2ac(q-1)}-\dfrac{\left(q^j-1\right)q^{-j}(c-aq)\left(acq^{j+1}-1\right)}{2ac\left(q^2-1\right)} \qquad &n=j, \\[1em]
        \dfrac{a+a^{-1}}{2}+\dfrac{(1-\alpha)(c-a)\left(q^{j+1}-1\right)q^{-j}\left(acq^j-1\right)}{2ac(q-1)}-\dfrac{\left(q^j-1\right)q^{-j}(c-aq)\left(acq^{j+1}-1\right)}{2ac\left(q^2-1\right)} \qquad &n=j+1, \\[1em]
       \end{cases} \label{bnodd}\\[1em] 
 u_n &= \begin{cases}
        \dfrac{\left(q^n-1\right) \left(q^n-q^{2 j+2}\right) \left(a c q^n-q\right) \left(q^n-a c q^{2 j+1}\right) \left(a q^n-c q^{j+1}\right) \left(c q^n-a q^{j+1}\right)}{4 a^2 c^2 \left(q^{j+1}+q^n\right)^2 \left(q^{2 n}-q^{2 j+1}\right) \left(q^{2 n}-q^{2 j+3}\right)} \qquad & n\neq j+1, \\[1em]
        \dfrac{(1-\alpha)\alpha(c-a)^2 q^{-2j}\left(q^{j+1}-1\right)^2 \left(acq^j-1\right)^2}{4a^2c^2(q-1)^2} \qquad &n=j+1.
       \end{cases} \label{unodd}
\end{align}
}%
{\em Remark :} It is straightforward to verify that these coefficients are persymmetric when $\alpha=1/2$, i.e. that the following relations hold : 
\begin{align} \label{persymmetry}
\begin{aligned}
 b_n &= b_{N-n}, &&n=0,1,\dots,N, \\
 u_n &= u_{N-n+1}, &&n=1,2,\dots,N.
\end{aligned}
\end{align}
For $\alpha\neq1/2$, the coefficients $b_j$ and $b_{j+1}$ are perturbed and no longer equal. 

Furthermore, owing to Favard's theorem, these polynomials will be orthogonal if they satisfy $u_n>0$ for $n=1,\dots,N$. With some easy computation, this is seen to be tantamount to the following conditions on the parameters 
\begin{align}
\begin{aligned}
  &0<q<1, \qquad 0<\alpha<1, \qquad c\neq a,\\
  &q<\tfrac{a}{c}<q^{-1}, \qquad ac<1 \quad \text{or} \quad ac>q^{1-N}.
\end{aligned}
\end{align}

\subsection{$q$-Difference equation}

Exploiting again the limit procedure given by \eqref{trunclimit}, it is possible to recover a $q$-difference equation for the $q$-para-Racah polynomials from that of the Askey-Wilson polynomials. Indeed, this procedure is trivial since the $q$-difference equation \eqref{AWdifferenceeqn} contains no parameters in the denominator. Thus, the $R_n(x)$ will satisfy
  \begin{align} \label{qdiffeqn}
   q^{-n}(1-q^n)(1-q^{n-N})R_n(x) = A(\theta)T_+R_n(x) -[A(\theta)+\overline{A}(\theta)]R_n(x) + \overline{A}(\theta)T_-R_n(x)
  \end{align}
  where the $T_\pm$ are again given by \eqref{qshiftop} and the coefficients $A(\theta)$ become
  \begin{align}
   A(\theta) = \frac{(1-a e^{i\theta})(1-a^{-1}q^{-j}e^{i\theta})(1-c e^{i\theta})(1-c^{-1}q^{-j} e^{i\theta})}{(1- e^{2i\theta})(1-q e^{2i\theta})}
  \end{align}
  with $\overline{A}(\theta)$ obtained by complex conjugation. The $q$-para-Racah polynomials are thus bispectral, but we remark that upon scaling the polynomials by a factor $q^{-n}$, the eigenvalues from equation \eqref{qdiffeqn} are degenerate in contrast to the usual classical orthogonal polynomials.

\subsection{Explicit expression}
It is possible to obtain an explicit expression for the $q$-para-Racah polynomials from the hypergeometric expression of the Askey-Wilson polynomials. Consider the series expansion \eqref{AWExplicit} and use the parametrization \eqref{trunclimit}. In the limit $t\to0$, the coefficients \eqref{AWExplicitCoefficients} reduces to 
\begin{align}
 \lim_{t\to0}A_{n,k} = \begin{cases}
                       \dfrac{(q^{-n},q^{n-2j-1};q)_k q^k}{(q^{-j},\tfrac{a}{c}q^{-j},ac,q;q)_k} \quad &\text{$k\le j$ and $k\le n$}, \\[1em]
                       \dfrac{(q^{n-2j-1};q)_{2j+1-n}(q;q)_{n+k-2j-2}(q^{-n};q)_k q^k}{\alpha (q^{-j};q)_j(q;q)_{k-j-1}(\tfrac{a}{c}q^{-j},ac,q;q)_k} \quad &\text{$k>j$ and $k\le n$}, \\[1em]
                       0 \quad &\text{otherwise}.
                       \end{cases}
\end{align}
This allows us to express the $R_n(x)$ as
\begin{align}
 R_n(x) = \eta_n \sum_k \left(\lim_{t\to0}A_{n,k}\right)\Phi_k(x)
\end{align}
where $\eta_n$ is a normalization factor to ensure the polynomials are monic. With the help of some well-known identities for $q$-Pochhammer symbols, the following expressions can be obtained : 
If $n<j$,
\begin{align}
R_n(x) = \eta_n \pGq{4}{3}{q^{-n},q^{n-2j-1},ae^{i\theta},ae^{-i\theta}}{q^{-j},ac,ac^{-1}q^{-j}}{q;q}. 
\end{align}
If $n=j$,
\begin{align}
R_j(x) = \eta_j \sum_{k=0}^j \frac{(q^{-j-1},ae^{i\theta},ae^{-i\theta};q)_k q^k}{(q,ac,ac^{-1}q^{-j};q)_k}. 
\end{align}
If $n=j+1$,
\begin{align}
R_{j+1}(x) = \eta_{j+1} \sum_{k=0}^j \frac{(q^{-j-1},ae^{i\theta},ae^{-i\theta};q)_k q^k}{(q,ac,ac^{-1}q^{-j};q)_k} + \eta_{j+1} \frac{(q^{-j-1},ae^{i\theta},ae^{-i\theta};q)_{j+1} q^{j+1}}{\alpha (q,ac,ac^{-1}q^{-j};q)_{j+1}}.
\end{align}
If $j+1<n\le N$,
{\small
\begin{align}
 R_n&(x) = \eta_n \pGq{4}{3}{q^{-n},q^{n-2j-1},ae^{i\theta},ae^{-i\theta}}{q^{-j},ac,ac^{-1}q^{-j}}{q;q} \\[1em] \notag
        &+ \eta_n \frac{(q^{n-2j-1};q)_{2j+1-n}(q^{-n},ae^{i\theta},ae^{-i\theta};q)_{j+1}(q;q)_{n-j-1}q^{j+1}}{\alpha (q^{-j};q)_j (q,ac,ac^{-1}q^{-j};q)_{j+1}}
        \pGq{4}{3}{q^{j+1-n},q^{n-j},aq^{j+1}e^{i\theta},aq^{j+1}e^{-i\theta}}{q^{j+2},acq^{j+1},ac^{-1}q}{q;q}. \label{RN} 
\end{align}
}%
The normalization $\eta_n$ is given by 
\begin{align}
 \eta_n = \begin{cases}
           \dfrac{(q,q^{-j},ac^{-1}q^{-j},ac;q)_n}{(q^{n-2j-1},q^{-n};q)_n(-2a)^n q^{n(n+1)/2}} \qquad &n\le j, \\[1em]
           \dfrac{\alpha(q^{-j};q)_j(q;q)_{n-j-1}(ac^{-1}q^{-j},ac,q;q)_n}{(q^{n-2j-1};q)_{2j+1-n}(q;q)_{2n-2j-2}(q^{-n};q)_n(-2a)^nq^{n(n+1)/2}} \qquad &n > j.
          \end{cases}
\end{align}
The $q$-para-Racah polynomials thus generally admit an explicit expression as a linear combination of two basic hypergeometric functions. However, due to some cancellations between parameters in the numerator and the denominator of the hypergeometric function when $n=j$ and $n=j+1$, the polynomials $R_j(x)$ and $R_{j+1}(x)$ have to be expressed as a sum which corresponds to a "truncated" hypergeometric series.  

Moreover, the non-monic $q$-para-Racah polynomials given by $R_n(x;a,c,\alpha|q)/\eta_n$ reduces to the (non-monic) para-Racah polynomials described in \cite{lemay2016racah} upon substituting
\begin{align}
 a\to q^a, \quad c\to q^c, \quad e^{i\theta}\to q^{ix}
\end{align}
and taking the limit $q\to1$. Akin to the connection between the Askey-Wilson and the Wilson polynomials, the $q\to1$ limit is taken in the explicit expressions of the polynomials instead of the recurrence relation.

\subsection{Orthogonality relation}
In order to obtain the orthogonality relation of the $q$-para-Racah polynomials, we begin by computing the characteristic polynomial which gives the spectrum of the Jacobi matrix or, equivalently, the orthogonality lattice. In a similar fashion to the derivation of the $R_n(x)$, consider    
\begin{align}
 \lim_{t\to0}(1-q^{(e_1+e_2)t})A_{N+1,k} = \begin{cases}
                                          \dfrac{(q^{-N-1};q)_k q^k}{\alpha (ac,ac^{-1}q^{-j};q)_k(q^{-j};q)_j(q;q)_{k-j-1}} \qquad &j+1\le k \le N+1, \\[1em]
                                          0 \qquad &\text{otherwise},
                                         \end{cases}
\end{align}
which can be summed with the $\Phi_k(x)$ to obtain 
\begin{align} \label{propR}
 R_{N+1}(x) \propto (ae^{i\theta},ae^{-i\theta};q)_{j+1} \times\pGq{3}{2}{q^{-j-1},aq^{j+1}e^{i\theta},aq^{j+1}e^{-i\theta}}{acq^{j+1},ac^{-1}q}{q;q}.
\end{align}
Note that it is possible to neglect the normalization constant since we are only interested in the zeros of $R_{N+1}(x)$. Now, using the Saalschutz $q$-summation formula, \eqref{propR} can be factorized as  
\begin{align} \label{CharacteristicP}
 R_{N+1}(x) \propto (ae^{i\theta},ae^{-i\theta};q)_{j+1}(ce^{i\theta},ce^{-i\theta};q)_{j+1}.
\end{align}
The orthogonality lattice will correspond to the zeros of \eqref{CharacteristicP} :   
\begin{align} \label{Olattice}
 x_{2s} &= \tfrac12(a^{-1}q^{-s}+aq^s) \qquad s=0,1,\dots,j, \\[1em]
 x_{2s+1} &= \tfrac12(c^{-1}q^{-s}+cq^s) \qquad s=0,1,\dots,j. 
\end{align}
Hence, the $q$-para-Racah polynomials will obey an orthogonality relation of the form
\begin{align}
 \sum_{s=0}^N w_s R_n(x_s)R_m(x_s) = u_1u_2\dots u_n \delta_{nm} 
\end{align}
where the $x_s$ are given by \eqref{Olattice} and the normalization constants are given by the recurrence coefficients \eqref{unodd}. A standard formula from the theory of orthogonal polynomials explicitly gives the weights \cite{1978_Chihara} :
\begin{align} \label{weightsformula}
 w_s = \frac{u_1u_2\dots u_N}{R_n(x_s)R'_{N+1}(x_s)}, \quad s=0,1,\dots,N.
\end{align}
However, due to the involved nature of $R_N(x)$ given by \eqref{RN}, we will use a simpler procedure which exploits the persymmetry that arises when $\alpha=1/2$ \cite{genest2017persymmetric}. In this case, when the polynomial $R_N(x)$ is evaluated at the zeros of the characteristic polynomial $x_s$, one obtains a simple expression which is due to the interlacing properties of their zeros :
\begin{align} \label{simpleRN}
 R_N(x_s) = \sqrt{u_1u_2\dots u_N}(-1)^{N+s}.
\end{align}
Combining \eqref{weightsformula} and \eqref{simpleRN}, it is easy to compute the weights for $\alpha=1/2$ which we shall denote by $\tilde{w}_s$. The weights for general $\alpha$ have been shown to be related to the $\tilde{w}_s$ by a simple multiplicative factor in \cite{genest2017persymmetric} :
\begin{align} \label{relatews}
 w_s \propto (1+\beta(-1)^s)\tilde{w}_s
\end{align}
where $\beta$ is a real parameter independent of $N$. One can easily obtain $\beta$ by comparing \eqref{weightsformula} and \eqref{relatews} for a fixed value of $N$ (e.g. $N=3$ for simplicity). Here, one obtains $\beta=1-2\alpha$. Carrying through the calculation, one readily obtains    
\begin{align}
 w_{2s} = -\frac{2(1-\alpha)K_N 2^{2j+1} a^j c^{j+1} q^{(2j+1)s+(j+1)j}\left(1-a^2 q^{2s}\right) \left(a^2;q\right)_{s} \left(q^{-j};q\right)_{s} (a c;q)_{s} \left(\frac{a q^{-j}}{c};q\right)_{s}}{(q;q)_j \left(a^2 q;q\right)_j \left(\frac{c}{a};q\right)_{j+1} (a c;q)_{j+1}\left(1-a^2\right) (q;q)_{s} \left(\frac{a q}{c};q\right)_{s} \left(a^2 q^{j+1};q\right)_{s} \left(a c q^{j+1};q\right)_{s}}
\end{align}
\begin{align}
 w_{2s+1} = \frac{2\alpha K_N 2^{2j+1} c^j a^{j+1} q^{(2j+1)s+(j+1)j}\left(1-c^2 q^{2s}\right) \left(c^2;q\right)_{s} \left(q^{-j};q\right)_{s} (ac;q)_{s} \left(\frac{c q^{-j}}{a};q\right)_{s}}{(q;q)_j \left(c^2 q;q\right)_j \left(\frac{a}{c};q\right)_{j+1} (ac;q)_{j+1}\left(1-c^2\right) (q;q)_{s} \left(\frac{cq}{a};q\right)_{s} \left(c^2 q^{j+1};q\right)_{s} \left(acq^{j+1};q\right)_{s}}
\end{align}
where $K_N$ is a normalization constant arising in the persymmetric case $\alpha=1/2$ and given by 
\begin{align}
 K_N = \sqrt{u_1u_2\dots u_N}.
\end{align}
It can easily be computed by using the persymmetry of the $u_n$ \eqref{persymmetry} :
\begin{align}
 K_N = \frac{(a-c) q^{-j} \left(q^{j+1}-1\right) \left(a c q^j-1\right) \left(q^{-j};q\right)_j^2 \left(q^{-2 j-1};q\right){}_j (q;q)_j (a c;q)_j \left(\frac{a q^{-j}}{c};q\right){}_j \left(\frac{c q^{-j}}{a};q\right){}_j \left(\frac{q^{-2 j}}{a c};q\right){}_j}{a c (q-1) 2^{2 j+2} \left(q^{-2 j-1};q^2\right){}_j \left(q^{-2 j};q^2\right)_j^2 \left(q^{1-2 j};q^2\right){}_j}.
\end{align}
The weights are normalized to verify  
\begin{align}
 \sum_{s=0}^j w_{2s} = 1-\alpha \qquad \sum_{s=0}^j w_{2s+1} = \alpha
\end{align}
which generalizes a known result for persymmetric Jacobi matrices. 

When $c=a$ the spectrum becomes doubly degenerate, i.e. $x_{2s}=x_{2s+1}$. This degeneracy is related with the degeneracy of the recurrence coefficient $u_{j+1}=0$ as seen from \eqref{unodd}. This means that the corresponding Hermitian Jacobi (tridiagonal) matrix of the recurrence coefficients $b_n$ and $u_n$ becomes reducible: it can be decomposed into a direct sum of two independent Jacobi matrices each having the same (now simple) spectrum $x_{2s}$.

\section{Even case : $N = 2j$}  

The construction of the $q$-para-Racah polynomials for even values of $N$ is similar to the $N$ odd case. We review quickly the procedure in this section and give the corresponding results. 

Let $N=2j$. The singular truncation \eqref{truncation} is achieved via the parametrization and limit
\begin{align} \label{trunclimit2}
 b = a^{-1}q^{-j+e_1t}, \qquad d = c^{-1}q^{-j+1+e_2t}, \qquad t\to0.
\end{align}
As before, the parameters $e_1$ and $e_2$ are not independent and can be encoded in a single deformation parameter $\alpha$ by
\begin{align}
 \frac{e_1}{e_1+e_2} = \alpha \qquad \frac{e_2}{e_1+e_2} = 1-\alpha.
\end{align}
\subsection{Recurrence relation}
The $q$-para-Racah polynomials $R_n(x;a,c,\alpha|q)$, or $R_n(x)$ for short, are defined by the recurrence relation
\begin{align} \label{recrel2}
 &x R_n(x) = R_{n+1}(x) + b_nR_n(x) + u_n R_{n-1}(x) 
\end{align}
with initial conditions $R_{-1}(x)=0, R_0(x)=1$ and with the coefficients given by 
\begin{align}
 &b_n = \lim_{t\to0}\tfrac{1}{2}(a+a^{-1}-A_n-C_n), \\
 &u_n = \lim_{t\to0}\tfrac{1}{4}A_{n-1}C_n. 
\end{align}
where $A_n$ and $C_n$ are the recurrence coefficients of the Askey-Wilson polynomials \eqref{rec_coef_wilson} in which we substituted parametrization \eqref{trunclimit2}. A straightforward calculation yields
\begin{align}
 b_n &= \dfrac{a+a^{-1}}{2}+\dfrac{\left(q^n-1\right) \left(a c q^{2 j}-q^n\right) \left(a q^{j+1}-c q^n\right)}{2ac \left(q^j+q^n\right) \left(q^{2 j+1}-q^{2 n}\right)}+\dfrac{\left(q^{2 j}-q^n\right) \left(a c q^n-1\right) \left(c q^j-a q^{n+1}\right)}{2ac \left(q^j+q^n\right) \left(q^{2 j}-q^{2 n+1}\right)}, \label{bneven}\\[1em]
 u_n &= \begin{cases}
          \dfrac{\left(q^n-1\right) \left(q^n-q^{2 j+1}\right) \left(a c q^n-q\right) \left(q^n-a c q^{2 j}\right) \left(a q^n-c q^j\right) \left(c q^n-a q^{j+1}\right)}{4 a^2 c^2 \left(q^j+q^n\right) \left(q^{j+1}+q^n\right) \left(q^{2 j+1}-q^{2 n}\right)^2} \quad & n\neq j,j+1, \\[1em]
          \dfrac{(1-\alpha) (c-a) q^{-2 j} \left(q^j-1\right) \left(q^{j+1}-1\right) (a q-c) \left(a c q^j-1\right) \left(a c q^j-q\right)}{4 a^2 c^2 (q-1)^2 (q+1)} \quad & n=j, \\[1em]
          \dfrac{\alpha     (c-a) q^{-2 j} \left(q^j-1\right) \left(q^{j+1}-1\right) (a q-c) \left(a c q^j-1\right) \left(a c q^j-q\right)}{4 a^2 c^2 (q-1)^2 (q+1)} \quad & n=j+1.
        \end{cases} \label{uneven}
\end{align}
The recurrence coefficients are also persymmetric, i.e. satisfy \eqref{persymmetry}, when $\alpha=1/2$. The positivity conditions $u_n>0$ for $n=1,\dots,N$ are verified when the parameters obey
\begin{align}
\begin{aligned}
  &0<q<1, \qquad 0<\alpha<1, \qquad c\neq a,\\
  &q<\tfrac{a}{c}<q^{-1}, \qquad ac<1 \quad \text{or} \quad ac>q^{1-N}.
\end{aligned}
\end{align}

\subsection{$q$-Difference equation}

The $q$-difference equation for $N=2j$ is obtained by inserting the limit procedure given by \eqref{trunclimit2} in \eqref{AWdifferenceeqn}. In this case, the $R_n(x)$ satisfy
  \begin{align} \label{qdiffeqn2}
   q^{-n}(1-q^n)(1-q^{n-N})R_n(x) = A(\theta)T_+R_n(x) -[A(\theta)+\overline{A}(\theta)]R_n(x) + \overline{A}(\theta)T_-R_n(x)
  \end{align}
  where the $T_\pm$ are given by \eqref{qshiftop} and the coefficients by
  \begin{align}
   A(\theta) = \frac{(1-a e^{i\theta})(1-a^{-1}q^{-j}e^{i\theta})(1-c e^{i\theta})(1-c^{-1}q^{-j+1} e^{i\theta})}{(1- e^{2i\theta})(1-q e^{2i\theta})}
  \end{align}
  and its complex conjugate. Again, the $q$-para-Racah polynomials are bispectral, but each eigenvalues is degenerate upon rescaling the polynomials by $q^{-n}$.
  
\subsection{Explicit expression}

An explicit expression for the $R_n(x)$ can readily be obtained by inserting the parametrization \eqref{trunclimit2} in the series expansion \eqref{AWExplicit}. Summing the resulting terms, one obtains for $n\le j$ :
\begin{align}
R_n(x) = \eta_n \pGq{4}{3}{q^{-n},q^{n-2j},ae^{i\theta},ae^{-i\theta}}{q^{-j},ac,ac^{-1}q^{-j+1}}{q;q}, 
\end{align}
and for $j+1\le n\le N$ :
\begin{align}
 R_n&(x) = \eta_n \pGq{4}{3}{q^{-n},q^{n-2j},ae^{i\theta},ae^{-i\theta}}{q^{-j},ac,ac^{-1}q^{-j+1}}{q;q} \\[1em] \notag
        &+ \eta_n \frac{(q^{n-2j};q)_{2j-n}(q^{-n},ae^{i\theta},ae^{-i\theta};q)_{j+1}(q;q)_{n-j}q^{j+1}}{\alpha (q^{-j};q)_j (q,ac,ac^{-1}q^{-j+1};q)_{j+1}}
        \pGq{4}{3}{q^{j+1-n},q^{n-j+1},aq^{j+1}e^{i\theta},aq^{j+1}e^{-i\theta}}{q^{j+2},acq^{j+1},ac^{-1}q^2}{q;q}
\end{align}
with monicity ensured by the normalization 
\begin{align}
 \eta_n = \begin{cases}
           \dfrac{(q,q^{-j},ac^{-1}q^{-j+1},ac;q)_n}{(q^{n-2j},q^{-n};q)_n(-2a)^n q^{n(n+1)/2}} \qquad &n\le j, \\[1em]
           \dfrac{\alpha(q^{-j};q)_j(q;q)_{n-j-1}(ac^{-1}q^{-j+1},ac,q;q)_n}{(q^{n-2j};q)_{2j-n}(q;q)_{2n-2j-1}(q^{-n};q)_n(-2a)^nq^{n(n+1)/2}} \qquad &n > j.
          \end{cases}
\end{align}
The polynomials of degree $j$ and $j+1$ need not be distinguished when $N=2j$ because the simplification of the parameters in the hypergeometric function does not change where the series truncate in constrast with the $N$ odd case.

\subsection{Orthogonality relation}

The characteristic polynomial can once more be computed via
\begin{align}
 R_{N+1}(x)= \eta_{N+1}\sum_{k=0}^\infty\lim_{t\to0}(1-q^{(e_1+e_2)t})A_{N+1,k}\Phi_k(x).
\end{align}
Carrying through the computation and using the Saalschutz $q$-summation formula, it can be expressed in factorized form as  
\begin{align}
 R_{N+1}(x) \propto (ae^{i\theta},ae^{-i\theta};q)_{j+1}(ce^{i\theta},ce^{-i\theta};q)_{j}.
\end{align}
The orthogonality grid is again has the form of a biexponential bi-lattice :
\begin{align}
 x_{2s} &= \tfrac12(a^{-1}q^{-s}+aq^s) \qquad s=0,1,\dots,j, \\
 x_{2s+1} &= \tfrac12(c^{-1}q^{-s}+cq^s) \qquad s=0,1,\dots,j-1, 
\end{align}
and the orthogonality relation is
\begin{align}
 \sum_{s=0}^N w_s R_n(x_s)R_m(x_s) = u_1u_2\dots u_n \delta_{nm}. 
\end{align}
As in the previous section, one can compute the general weights by using the persymmetry when $\alpha=1/2$. The result is 
\begin{align}
 w_{2s} = \frac{(1-\alpha)K_N a^j c^{j} q^{(2j)s}\left(1-a^2 q^{2s}\right) \left(a^2;q\right)_{s} \left(q^{-j};q\right)_{s} (a c;q)_{s} \left(\frac{a q^{-j+1}}{c};q\right)_{s}}{(q;q)_{j} \left(a^2 q;q\right)_j \left(\frac{c}{a};q\right)_{j+1} (a c;q)_{j}\left(1-a^2\right) (q;q)_{s} \left(\frac{a q}{c};q\right)_{s} \left(a^2 q^{j+1};q\right)_{s} \left(a c q^{j};q\right)_{s}},
\end{align}
\begin{align}
 w_{2s+1} = \frac{-\alpha K_N  a^{j+1}c^{j-1}  q^{(2j)s}\left(1-c^2 q^{2s}\right) \left(c^2;q\right)_{s} \left(q^{-j+1};q\right)_{s} (ac;q)_{s} \left(\frac{c q^{-j}}{a};q\right)_{s}}{(q;q)_{j-1} \left(c^2 q;q\right)_{j-1} \left(\frac{a}{c};q\right)_{j+1} (ac;q)_{j+1}\left(1-c^2\right) (q;q)_{s} \left(\frac{cq}{a};q\right)_{s} \left(c^2 q^{j};q\right)_{s} \left(acq^{j+1};q\right)_{s}},
\end{align}
with
{\small
\begin{align}
 K_N = \frac{q^{2j^2+1} (c-a) (1+q^j) (1-q^{j+1}) (1-q^{2 j+1}) (c-a q) (q;q)_j (q^{-2 j-1};q)_j (\frac{a c}{q};q)_{j+2} (\frac{c q^{-j-1}}{a};q)_j (\frac{q^{-2j}}{ac};q)_j (\frac{a q^{-j}}{c};q)_j}{(1-q)^2 (a c-q) (q^{-2 j-1};q^2)_j^2 (-q;q)_j^2 (a-c q^j) (1-a c q^{2j}) (c-a q^{j+1})}.
\end{align}
}%
The weights also satisfy the relations
\begin{align}
 \sum_{s=0}^j w_{2s} = 1-\alpha, \qquad \sum_{s=0}^{j-1} w_{2s+1} = \alpha.
\end{align}

\section{Special cases}
\subsection{$q$-para-Krawtchouk}
Under an appropriate reparametrization, it is possible to reduce the $R_n(x)$ to polynomials orthogonal with respect to an exponential bi-lattice instead of a biexponential bi-lattice. We call the corresponding polynomials the $q$-para-Krawtchouk polynomials because they reduce to the para-Krawtchouk polynomials when $q\to1$. The persymmetric case $\alpha=1/2$ for $N=2j+1$ has been briefly mentionned in \cite{tsujimoto2017tridiagonal}. We here obtain the general $q$-para-Krawtchouk polynomials for general $N$ and general $\alpha$. To this end, let us first rewrite the parameters $a$ and $c$ in terms of new parameters $\theta$ and $\Delta$ given by
\begin{align}
 \theta = ac, \qquad \Delta = \frac{a}{c}, \\
 a^2 = \theta\Delta, \qquad c^2 = \frac{\theta}{\Delta}.
\end{align}
Rescaling the lattice $x_s$ and taking the limit $\theta\to\infty$, one obtains an exponential bi-lattice $y_s$ described by  
\begin{align}
 y_{2s}= \lim_{\theta\to\infty} \frac{2a}{\theta} x_{2s} = \Delta q^s, \\
 y_{2s+1}= \lim_{\theta\to\infty} \frac{2a}{\theta} x_{2s+1} =  q^s. 
 \end{align}
The $q$-para-Krawtchouk polynomials $Q_n(y)$ can be obtained by taking a similarity transformation
\begin{align}
 x = \frac{\theta}{2a} y, \qquad Q_n(y) = \left(\frac{\theta}{2a}\right)^{-n} R_n(x).
\end{align}
Under this transformation, the recurrence relations \eqref{recrel1} and \eqref{recrel2} of the $q$-para-Racah becomes
\begin{align}
 y Q_n(y) = Q_{n+1} + \tilde{b}_n Q_n(y) + \tilde{u}_n Q_{n-1}(y)
\end{align}
where
\begin{align}
 \tilde{b}_n = \lim_{\theta\to\infty} \frac{2a}{\theta} b_n, \\
 \tilde{u}_n = \lim_{\theta\to\infty} \frac{4a^2}{\theta^2} u_n. 
 \end{align}
Inserting the recurrence coefficients \eqref{bnodd} and \eqref{unodd} in the previous formula gives
\begin{align}
 \tilde{b}_n = \begin{cases}
                \frac{q^{n+j}(1+q^{j+1})(1+\Delta)}{(q^j+q^n)(q^{j+1}+q^n)} \quad &n\neq j,j+1, \\[0.5em]
                \Delta - \frac{\alpha(1-q^{j+1})(\Delta-1)}{1-q} +\frac{q(1-q^j)(\Delta q-1)}{(1-q^2)} \quad &n=j, \\[0.5em]
                \Delta - \frac{(1-\alpha)(1-q^{j+1})(\Delta-1)}{1-q} +\frac{q(1-q^j)(\Delta q-1)}{(1-q^2)} \quad &n=j+1, 
               \end{cases}\\[1em]       
 \tilde{u}_n = \begin{cases}
                \frac{q^{2j+1+n}(1-q^n)(q^{2j+2}-q^n)(q^n-\Delta q^{j+1})(q^{j+1}-\Delta q^{n})}{(q^{j+1}+q^n)^2(q^{2j+1}-q^{2n})(q^{2j+3}-q^{2n})} \quad &n\neq j+1, \\[0.5em]
                \frac{\alpha(1-\alpha)(\Delta-1)^2(1-q^{j+1})^2}{(1-q)^2}\quad &n=j+1 
               \end{cases}               
\end{align}
for $N=2j+1$. Using instead \eqref{bneven} and \eqref{uneven} gives
\begin{align}
 \tilde{b}_n = \Delta - \frac{q^{2j}(q^n-1)(q^n-\Delta q^{j+1})}{(q^j+q^n)(q^{2j+1}-q^{2n})} + \frac{q^{n}(q^{2j}-q^n)(q^j-\Delta q^{n+1})}{(q^j+q^n)(q^{2j}-q^{2n+1})},\\[1em]
 \tilde{u}_n = \begin{cases}
                \frac{q^{2j+n}(q^n-1)(q^{2j+1}-q^n)(q^n-\Delta q^{j+1})(\Delta q^n-q^j)}{(q^j+q^n)(q^{j+1}+q^n)(q^{2j+1}-q^{2n})^2} \quad &n\neq j,j+1, \\[0.5em]
                \frac{(1-\alpha)(1-q^j)(1-q^{j+1})(\Delta-1)(1-q\Delta)}{(1-q)^2(1+q)} \quad &n=j, \\[0.5em]
                \frac{\alpha(1-q^j)(1-q^{j+1})(\Delta-1)(1-q\Delta)}{(1-q)^2(1+q)} \quad &n=j+1
               \end{cases}
\end{align}
for $N=2j$. In the limit $q\to1$, this recurrence relation reduces to the one of the para-Krawtchouk polynomials up to another similarity transformation. Remark that these coefficients are also persymmetric when $\alpha=1/2$. The orthogonality relation for the $Q_n(y)$ can be obtained by the same procedure used in section 3 and 4 to obtain the orthogonality relation for the $q$-para-Racah polynomials. Omitting the details, one readily finds 
\begin{align}
 \sum_{s=0}^N w_s Q_n(y_s)Q_m(y_s) = \tilde{u}_1\tilde{u}_2\dots\tilde{u}_n\delta_{nm}
\end{align}
where the weights are given by 
\begin{align}
 w_{2s}&=K_N\frac{(1-\alpha) \left(1-\frac{1}{\Delta }\right) q^s \left(\Delta  q^{-j};q\right)_j \left(\frac{q^{-j}}{\Delta };q\right)_j \left(q^{-j};q\right)_s \left(\Delta  q^{-j};q\right)_s}{(q;q)_s \left(\frac{1}{\Delta };q\right)_{j+1}\Delta ^{j} (\Delta  q;q)_s} \\
 w_{2s+1}&=K_N\frac{\alpha  (1-\Delta ) \Delta ^j q^s \left(\frac{q^{-j}}{\Delta };q\right)_j \left(\Delta  q^{-j};q\right)_j \left(q^{-j};q\right)_s \left(\frac{q^{-j}}{\Delta };q\right)_s}{(q;q)_s (\Delta ;q)_{j+1} \left(\frac{q}{\Delta };q\right)_s}\\
 K_N &= \frac{(-1)^j q^{j (j-1)} \left(1-q^{2 j+1}\right)}{(1-q) (-q;q)_j \left(q^{-2 j-1};q^2\right)_j}.
\end{align}
for $N=2j+1$ and by 
\begin{align}
 w_{2s} &= K_N\frac{(1-\alpha ) q^s \left(q^{-j};q\right)_s \left(\Delta  q^{1-j};q\right)_s}{\Delta ^j (q;q)_s \left(\frac{q}{\Delta };q\right)_j (\Delta  q;q)_s} \\
 w_{2s+1} &= K_N\frac{\alpha  \Delta ^{j-1} \left(1-q^j\right) q^s \left(q^{1-j};q\right)_s \left(\frac{q^{-j}}{\Delta };q\right)_s}{\left(1-\frac{q^j}{\Delta }\right) (q;q)_s (\Delta  q;q)_j \left(\frac{q}{\Delta };q\right)_s} \\
 K_N &= \frac{(-1)^j q^{\frac{3}{2} j (j-1)} \left(q^j+1\right) \left(1-q^{j+1}\right) \left(1-q^{2 j+1}\right) \left(q^{-2 j-1};q\right)_j (1-\Delta  q) \left(\frac{q^{-j-1}}{\Delta };q\right)_j \left(\Delta  q^{-j};q\right)_j}{\left(1-\Delta  q^{j+1}\right) (1-q)^2 \left(q^{-2 j-1};q^2\right)_j^2 (-q;q)_j^2}
\end{align}
for $N=2j$. 

\subsection{Reduction to a single lattice}
Consider the persymmetric case $\alpha=\frac12$. For $c = a q^{\frac12}$, the orthogonality lattice reduces to a single biexponential lattice of the form 
\begin{align}
 x_s = \frac12 (a^{-1}q^{-\frac{s}{2}}+a q^{\frac{s}{2}}) \quad s=0,1,2,\dots,N.
\end{align}
In this setting, the $q$-para-Racah polynomials connect with the $q$-Racah polynomials in base $q^{\frac12}$. More precisely, the following relation holds for any $N$ : 
\begin{align} \label{qRacah}
R_n(x;a,aq^{\frac12},\tfrac12|q) =(2a)^{-n} p_n(2ax;aq^{-\frac14},-a^{-1}q^{-j-\frac34},-aq^{-\frac14},-aq^{-\frac14}|q^{\frac12})  
\end{align}
where $p_n(y;\alpha,\beta,\gamma,\delta|q) \equiv p_n(y)$ are the monic $q$-Racah polynomials defined in \cite{2010_Koekoek_&Lesky&Swarttouw}. This can be checked by substituting directly \eqref{qRacah} in the recurrence relation \eqref{recrel1} for $N$ odd and \eqref{recrel2} for $N$ even and comparing the coefficients with those of the monic $q$-Racah polynomials found in \cite{2010_Koekoek_&Lesky&Swarttouw}. 

In addition, this special instance of $q$-para-Racah polynomials bears a connection with the dual-Hahn polynomials in the limit $q\to1$. To see this, consider the recurrence coefficients $A_n$ and $C_n$ given in \eqref{reccoefA} and \eqref{reccoefC} and let $\alpha=\frac12$ and $c=aq^{\frac12}$. Now, substitute $a\to q^a$ and compute the limits 
\begin{align}
\begin{aligned}
 \lim_{q\to1} \frac{A_n}{(1-q^{\frac12})^2} &= (n+\tfrac{4a-1}{2}+1)(n-N), \\
 \lim_{q\to1} \frac{C_n}{(1-q^{\frac12})^2} &= n(n-\tfrac{4a-1}{2}-N-1).
\end{aligned}
\end{align}
These results corresponds precisely to the recurrence coefficients (also denoted by $A_n$ and $C_n$) of the dual-Hahn polynomials given in \cite{2010_Koekoek_&Lesky&Swarttouw} with parameters $\gamma = \delta = \frac{4a-1}{2}$. It is a trivial matter to verify that the same result holds for even values of $N$. This is in perfect correspondance with the special case of the para-Racah polynomials that reduces to the dual-Hahn polynomials with parameters $\gamma = \delta = \frac{4a-1}{2}$ when the orthogonality bi-lattice of the former reduces to a single lattice \cite{lemay2016racah}.       

\section{Conclusion}

To summarize, we constructed new orthogonal polynomials from a singular truncation of the Askey-Wilson polynomials. They have been called the $q$-para-Racah polynomials because their construction is parallel to the one of the para-Racah polynomials starting from the Wilson polynomials. Furthermore, they can be connected by a $q\to1$ limit in their (unnormalized) explicit expression in a similar fashion to the connection between the Askey-Wilson and the Wilson polynomials. A three-term recurrence relation, a $q$-difference equation, an explicit expression and the orthogonality relation have been obtained both for sets containing an even or odd numbers of polynomials. We further characterized the $q$-para-Krawtchouk polynomials as a special case of the $q$-para-Racah polynomials. This is also of interest because these last polynomials had never been much characterized in the literature before. This specialization occurs in a limit where the orthogonality grid reduces from a biexponential bi-lattice to an exponential bi-lattice. A connection to the $q$-Racah and dual-Hahn polynomials has also been presented in the special case where the bi-lattice reduces to a single lattice. 

The $q$-para-Racah and the $q$-para-Krawtchouk are both associated to an isospectral deformation of persymmetric Jacobi matrices. Specifically, their Jacobi matrices are persymmetric only when the parameter $\alpha=1/2$. This is an interesting feature which we hope could see them arise in future applications. An idea is to interpret their Jacobi matrices as the restriction to the one-excitation sector of an Heisenberg spin chain Hamiltonian and to study their ability to produce transfer of quantum information or generate entangled pairs. Another direction would be to study their bispectrality. Although the $q$-para-Racah polynomials possess the bispectrality property, the spectrum of the corresponding $q$-difference operators in \eqref{qdiffeqn} and \eqref{qdiffeqn2} is doubly degenerate. This means that the $q$-para-Racah polynomials do not belong to the category of "classical" orthogonal polynomials with the Leonard duality property. It would be interesting to find an appropriate algebraic description of these polynomials. We hope to report on these questions in the near future.

\section*{Acknowledgments}
JML holds an Alexander-Graham-Bell PhD fellowship from the Natural Science and Engineering Research Council (NSERC) of Canada. LV is grateful to NSERC for support through a discovery grant.

\bibliographystyle{elsarticle-num}
\bibliography{Bibliography_ParaRacah.bib}

\end{document}